# Semi-Selfdecomposable Laws and Related Processes

S Satheesh

NEELOLPALAM, S. N. Park Road
Trichur – 680 004, **India.**

ssatheesh@sancharnet.in

and

E Sandhya

Department of Statistics, Prajyoti Niketan College
Pudukkad, Trichur – 680 301, **India.**

esandhya@hotmail.com

**Abstract.** In this note we identify the class of distributions for $X_n$ that can generate a linear, additive, first order auto-regressive scheme $X_n = \rho X_{n-1} + \varepsilon_n$, $0<\rho<1$, where $\{X_n\}$ is marginally stationary as semi-selfdecomposable laws. We give a method to construct these distributions. Its implications in subordination and selfdecomposability of Levy processes are given. The discrete analogues of these processes are also discussed.

## 1. Introduction

This investigation is motivated by the invited talk by Dr. N Balakrishna of CUSAT, at the annual conference of ISPS. Discussing non-Gaussian time series models Balakrishna (2004) mentioned the connection of selfdecomposable (SD) laws to the first order auto-regressive (AR(1)) model referring to Gaver and Lewis (1980). The linear additive AR(1) scheme considered here is described by the random variables (r.v) $\{X_n\}$ and $\{\varepsilon_n\}$ where $\{\varepsilon_n\}$ are independent and identically distributed (i.i.d), satisfying

$$X_n = \rho X_{n-1} + \varepsilon_n, \text{ some } 0<\rho<1. \tag{1}$$

Assuming $\{X_n\}$ to be marginally stationary, (1) can be written in terms of their characteristic functions (CF) as:

$$f_X(s) = f_X(\rho s) f_\varepsilon(s), \text{ for some } 0<\rho<1. \tag{2}$$

Now, a CF $f$ is SD if for every $c \in (0,1)$ there exists a CF $f_c$ such that

$$f(s) = f(cs) f_c(s), \forall s \in \mathbf{R.} \tag{2a}$$

Gaver and Lewis (1980, p.736) did record that this requirement (*ie.* (2a) is to be satisfied <u>for every</u> $c \in (0,1)$) in the theory of SD r.vs is its limitation in its relation to the stationary AR(1) model, and though the full range is desirable it may not occur.

Thus SD laws can generate the stationary AR(1) sequence $\{X_n\}$ but the distributions of such an $\{X_n\}$ need not necessarily be SD. Then how large is this class? Can we construct such distributions? Answer to these questions also gives us the possibility of constructing a wide range of AR(1) models that are non-Gaussian. An offshoot of the method of construction presented here is certain results on the subordination of processes with additive increments. These are the problems that we discuss here.

To sketch the background, Maejima and Naito (1998) conceived semi-SD laws on $\mathbf{R}^d$, $d>0$ integer, by considering the limit laws of subsequences of normalized partial sums of independent r.vs satisfying the infinitesimal condition and where the subsequence $\{k_n\}$ satisfies $k_n/k_{n+1} \to c<1$. In this note we restrict our discussion to **R**.

**Definition.1** (Maejima and Naito, 1998). A probability distribution $\mu$ on **R**, with CF $f$ is semi-SD($b$) if for some $b\in(0,1)$ there exists a CF $f_o$ such that

$$f(s) = f(bs) f_o(s), \ \forall \ s\in \mathbf{R}. \tag{3}$$

**Definition.2** (Pillai, 1971; Maejima and Sato, 1999). A CF $f$ that is infinitely divisible (ID) is semi-stable($a,b$) (this notion is by Levy) if for all $s\in \mathbf{R}$ and for some $0<b<1<a$,

$$f(u) = exp\{-\psi(u)\}, \text{ where } \psi(u) = a\psi(bu). \tag{4}$$

We do not consider the general case $0<|b|<1$ here. It is known (Pillai and Anil (1997) that here $\psi(u) = |u|^\alpha h(u)$, $\alpha \in (0,2]$, where $h(u)$ is a periodic function in $ln(u)$ with period $-ln(b)$ and $ab^\alpha =1$. When $h(u)$ is a constant $f$ is stable.

A stochastically continuous process $\{X(t), t\geq 0\}$ having stationary independent and additive increments and $X(o) = 0$ is called a Levy process. A Levy process $X(t)$ such that $X(1)$ is SD (semi-SD) will be called a SD (semi-SD) process. Since SD and semi-SD are ID the above Levy processes are well defined. Selfsimilar (SS) processes are processes that are invariant in distribution under suitable scaling of time and space. From Maejima and Sato (1999) and Embrechts and Maejima (2000) we have: A process $\{X(t), t\geq 0\}$ is SS if for any $b>0$ there is a unique exponent $H>0$ such that $\{X(bt)\} \stackrel{d}{=} \{b^H X(t)\}$. If the above relation holds for some $b>0$ only, then $\{X(t)\}$ is semi-SS. A Levy process $\{X(t)\}$ is SS (semi-SS) iff the distribution of $X(1)$ is stable (semi-stable). We will refer to these Levy processes as stable and semi-stable processes.

**Note.1** Here we consider the notions of stability, semi-stability, selfsimilarity and semi-selfsimilarity in the strict sense only.

**Note.2** Maejima and Sato (1999) also showed that for any given $H>0$ and a SD (semi-SD) law there exists a SS (semi-SS) process $\{X(t)\}$ with independent increments such that the distribution of $X(1)$ is SD (semi-SD).

**Remark.1** Without loss of generality we may consider the range of $b$ as $0<b<1$ in the description of selfsimilarity because it is equivalent to $\{(b^{-1})^H X(t)\} \stackrel{d}{=} \{X(b^{-1}t)\}$ and thus the whole range of $b>0$ is covered. Also the range $0<b<1$ appears more meaningful in its interpretation as a scaling factor and in its role in semi-stable, SD and semi-SD processes.

Subordination of processes is discussed in Feller (1971, p.573). Let $\{Y(t), t \geq 0\}$ be a Levy process with CF $exp\{-t\psi(u)\}$ and $\{T(t), t \geq 0\}$ a positive Levy process independent of $\{Y(t)\}$ with Laplace transform (LT) $\varphi^t$. Then the process $\{X(t), t \geq 0\}$ is said to be subordinated to $\{Y(t)\}$ by the directing process $\{T(t)\}$ if $\{X(t)\} \stackrel{d}{=} \{Y(T(t))\}$ and the CF of $\{X(t)\}$ is given by

$$f = \{\varphi(\psi)\}^t. \tag{5}$$

Here we are randomizing the time parameter $t$ of $\{Y(t)\}$ by $\{T(t)\}$.

Satheesh (2002 and 2004) discussed φ-ID, φ-semi-stable and φ-stable laws as follows. If $\varphi$ is a LT then the CF $\varphi\{-ln\ \omega(s)\}$ is φ-ID (φ-semi-stable and φ-stable) if the CF $\omega(s)$ is ID (semi-stable and stable). These are also the φ-mixture of the CF $\omega(s)$ that is ID (semi-stable and stable). Sandhya (1991) has discussed the connection between subordination and geometrically ID laws. Satheesh (2002) showed: For every $\varphi$-mixture of an ID law with CF $h(s) = \varphi\{-\ln \omega(s)\}$, the CF $h^t$ is that of a stochastic process $\{Y(T(t)), t \geq 0\}$ that is subordinated to the process $\{Y(t), t \geq 0\}$ having stationary and independent increments with CF $\omega^t$, by the directing process $\{T(t), t \geq 0\}$ with LT $\varphi^t$. Sato (2001) showed that: If $\{Y(t), t \geq 0\}$ is a stable (strictly) process and $\{T(t), t \geq 0\}$ is a SD process then the process $\{X(t), t \geq 0\}$ subordinated to $\{Y(t)\}$ directed by $\{T(t)\}$ is also SD. This is the process version of the method for constructing SD laws in Satheesh (2002, 2004) stated here as theorem.2, the proof of which is simpler than that in Sato (2001). Other works in this context are Shanbhag and Sreehari (1979) and Iksanov and Jurek (2003).

With this background this note is organized as follows. In the next section we prove that $\{X_n\}$ defines an AR(1) sequence iff $X_n$ is semi-SD and then show how to construct semi-SD laws. Implications of these constructions in the context of stochastic processes are then discussed generalizing and extending the results of Sato (2001) and Satheesh (2002). We then describe the discrete case also.

## 2. Results

**Theorem.1** A sequence $\{X_n\}$ of r.vs defines an AR(1) sequence that is marginally stationary with $0<\rho<1$ iff $X_n$ is semi-SD($\rho$).

*Proof.* Suppose $X_n$ is semi-SD($\rho$) with CF $f$. Then there exists a r.v $\varepsilon_n$ with CF $f_o$ such that the CFs of $X_n$ and $\varepsilon_n$ are related by

$f(s) = f(\rho s) f_o(s)$, for some $0<\rho<1$.

Thus from (2) $\{X_n\}$ defines a marginally stationary AR(1) sequence with $0<\rho<1$, and innovation sequence $\{\varepsilon_n\}$ with CF $f_o$. The converse is clear by comparing (2) and (3).

**Theorem.2** φ-(strictly) stable laws are SD if $\varphi$ is SD.

This was proved in Satheesh (2002, 2004). Generalizing this we have:

**Theorem.3** φ-semi-stable($a,b$) laws are semi-SD($b$) if $\varphi$ is SD.

*Proof.* If $\varphi$ is SD then the CFs of φ-semi-stable($a,b$) laws can be written as:

$$\varphi\{\psi(s)\} = \varphi\{c\psi(s)\}\varphi_o\{\psi(s)\}, \text{ for every } 0<c<1 \text{ and for all } s\in \mathbf{R},$$

where $\psi(s) = a\psi(bs)$, for some $0<b<1<a$, and $0<\alpha\le 2$ satisfying $ab^\alpha = 1$. (6)

Hence when $c = b^\alpha$ we have:

$$\varphi\{\psi(s)\} = \varphi\{\psi(bs)\}\varphi_o\{\psi(s)\} \text{ for some } 0<b<1, \text{ completing the proof.}$$

**Theorem.4** φ-semi-stable($a,b$) laws are semi-SD($b$) if $\varphi$ is semi-SD($b^\alpha$), $0<\alpha\le 2$. The Proof follows by a similar line of argument as above.

**Example.1** Consider the generalized semi-$\alpha$-Laplace law with CF $\{1+\psi(s)\}^{-\beta}$, where $\psi(s) = a\psi(bs)$ as in (6) and $\beta>0$. Since the gamma(1,$\beta$) law is SD and the above CF is a gamma mixture of semi-stable($a,b$) laws the generalized semi-$\alpha$-Laplace law is semi-SD($b$). Now write $b=v^\alpha$ for $v\in(0,1)$ and $\alpha\in(0,2]$. Setting $\varphi$ to be this semi-SD($v^\alpha$) law the corresponding φ-semi-stable($a,v$) laws are semi-SD($v$).

From theorems 2, 3 and 4 we now have the following results (their process versions) generalizing and extending the results of Sato (2001) and Satheesh (2002) in the context of subordination of stochastic processes. These results follow from the equation (5) and respectively from theorems 2, 3 and 4.

**Theorem.5** The stochastic process subordinated to a (strictly) stable process by a directing SD process is SD (Sato, 2001).

**Theorem.6** The stochastic process subordinated to a semi-stable($a$,$b$) process by a directing SD process is semi-SD($b$).

**Theorem.7** The stochastic process subordinated to a semi-stable($a$,$b$) process by a directing semi-SD($b^{\alpha}$) process is semi-SD($b$).

We know that stable laws are SD. We now prove that;

**Theorem.8** Semi-stable($a$,$b$) laws are semi-SD($b$).

*Proof.* From definition.2 we have: A CF $f(s)$ is semi-stable if it is ID and satisfies

$f(s) = \{f(bs)\}^a$, for some $0<b<1<a$ and $\forall\ s \in \mathbf{R}$. That is,

$f(s) = f(bs)\{f(bs)\}^{a-1} = f(bs)\,f_o(s)$.

That proves the assertion by Definition.1.

To have analogous results in the discrete setup we need the following lemma.

**Lemma.1** (Satheesh and Nair, 2002) If $\phi(s)$, $s>0$ is a LT, then $P(s) = \phi(1-s)$, $0<s<1$ is a probability generating function (PGF). Conversely, if $P(s)$ is a PGF and $P(1-s)$ is completely monotone for all $s>0$, then $\phi(s) = P(1-s)$ is a LT.

Using this lemma next we describe the PGF of a discrete semi-SD law and justify it by deriving it from the LT of a semi-SD law.

**Theorem.9** An integer-valued distribution with PGF $P(s)$ is discrete semi-SD($c$) if for some $0<c<1$, there exists another PGF $P_o(s)$ such that

$$P(s) = P(1-c+cs)\, P_o(s),\ \forall\ s\in(0,1).$$

*Proof.* Semi-SD(c) laws on $[0,\infty)$ are defined by LTs $\phi(s)$, satisfying,

$$\phi(s) = \phi(cs)\,\phi_o(s),\ \forall\ s>0 \text{ and some } 0<c<1$$

where $\phi_o(s)$ is another LT. In terms of PGFs (constructed by Lemma.1 from these LTs) this equation reads, for some $0<c<1$,

$$P(1-s) = P(1-cs)\, P_o(1-s),\ \forall\ s\in(0,1).$$

Setting $1-s = u$ in this equation we reach the conclusion.

**Corollary.1** If a LT $\phi(s)$ is semi-SD($c$) then the PGF $P(s) = \phi(1-s)$ is discrete semi-SD($c$). Conversely if a PGF $P(s)$ is discrete semi-SD($c$) and $\phi(s) = P(1-s)$ is a LT, then $\phi(s)$ is semi-SD($c$).

**Theorem.10** Discrete semi-stable($a,b$) laws are discrete semi-SD($b$).

Al-Osh and Alzaid (1987) had considered integer-valued AR(1) sequences. As an integer-valued analogue of theorem.1 we can now characterize (using theorem.9) integer-valued AR(1) sequences that are marginally stationary as:

**Theorem.11** A sequence $\{X_n\}$ of integer-valued r.vs defines a marginally stationary AR(1) sequence with $0<\rho<1$ iff $X_n$ is discrete semi-SD($\rho$).

Using lemma.1 discrete semi-stable and discrete stable laws have been discussed in Satheesh and Nair (2002). The discrete analogues of φ-ID, φ-semi-stable and φ-stable laws are described as follows. If $\varphi$ is a LT, then the PGF $\varphi\{-ln\ Q(s)\}$ is φ-ID (φ-semi-stable, φ-stable) if the PGF $Q(s)$ is ID (semi-stable, and stable). Now results analogous to theorems 2, 3, 4, 5, 6 and 7 can be stated and proved. We just state the analogues of theorems 4 and 7.

**Theorem.12** φ-semi-stable($a$,$b$) PGFs are discrete semi-SD($b$) if $\varphi$ is semi-SD($b^\alpha$), $0<\alpha\leq 2$. Consequently, the discrete state stochastic process subordinated to a discrete state semi-stable($a$,$b$) process by a directing semi-SD($b^\alpha$) process is discrete semi-SD($b$).

Thus we can construct a variety of distributions useful in AR(1) modeling, SS and semi-SS processes.

**References.**

**Al-Osh, M. A. and Alzaid, A. A. (1987).** First order integer-valued autoregressive (INAR(I)) process, *J. Time Ser. Anal.*, **VIII**, 261-275.

**Balakrishna, N. (2004).** Estimation in non-Gaussian time series, invited talk at the 24[th] Annual Conference of the *Indian Society for Probability and Statistics*, held at Pala, Kerala, India during 4-6 November 2004.

**Embrechts, C. and Maejima, M. (2000).** An Introduction to the theory of selfsimilar stochastic processes, *Int. J. Mod. Phy.* **B14**, 1399-1420.

**Feller, W. (1971).** *An Introduction to Probability Theory and Its Applications*, Wiley, New York, Vol.**2**, 2[nd] Edition.

**Gaver, D. P. and Lewis, P. A. W. (1980).** First-order autoregressive gamma sequences and point processes, *Adv. Appl. Probab.*, **12**, 727-745.

**Iksanov, A. M. and Jurek, Z. J. (2003).** Shot noise distributions and selfdecomposability, *Stoch. Anal. Appl.*, **21**, 593 – 609.

**Maejima, M. and Naito, Y. (1998).** Semi-selfdecomposable distributions and a new class of limit theorems, *Probab. Theor. Rel. Fields*, **112**, 13-31.

**Maejima, M. and Sato, K. (1999).** Semi-selfsimilar processes, *J. Theor. Probab.*, **12**, 347-383.

**Pillai, R. N. (1971).** Semi-stable laws as limit distributions, *Ann. Math. Statist.*, **42**, 780–783.

**Pillai, R. N. and Anil, V. (1996).** Symmetric stable, $\alpha$-Laplace, Mittag-Leffler and related laws/ processes and the integrated Cauchy functional equation, **34**, 97-103.

**Sandhya, E. (1991).** On geometric infinite divisibility, *p*-thinning and Cox processes, *J. Kerala Statist. Assoc.*, **7**, 1-10.

**Satheesh, S. (2002).** Aspects of randomization in infinitely divisible and max-infinitely divisible laws, *Probstat Models*, **1**, June-2002, 7-16.

**Satheesh, S. (2004).** Another look at random infinite divisibility, *Statist. Meth.*, **6**, 123 –144.

**Satheesh, S. and Nair, N. U. (2002).** Some classes of distributions on the non-negative lattice. *J. Ind. Statist. Assoc.*, **40**, 41–58.

**Sato, K. (2001).** Subordination and selfdecomposability, *Statist. Probab. Lett.*, **54**, 317 – 324.

**Shanbhag, D. N. and Sreehari, M. (1979).** An extention of Goldie's result and further results in infinite divisibility, *Z. Wahr. Verw. Geb.*, **47**, 19 – 25.